\documentclass[aps,pre,onecolumn,showpacs,amsfonts,amsmath,amssymb,superscriptaddress,nofootinbib,12pt]{revtex4-1}

\usepackage[dvips]{graphicx}	
 
\usepackage{txfonts}
\usepackage{hyperref}
\usepackage[usenames,dvipsnames]{color}
\usepackage{subfigure}

\usepackage{float}
\usepackage{empheq} 

\graphicspath{.,{figures/}}

\newcommand{\BEAS}{\begin{eqnarray*}}
\newcommand{\EEAS}{\end{eqnarray*}}
\newcommand{\BEA}{\begin{eqnarray}}
\newcommand{\EEA}{\end{eqnarray}}
\newcommand{\BEQ}{\begin{equation}}
\newcommand{\EEQ}{\end{equation}}
\newcommand{\bit}{\begin{itemize}}
\newcommand{\eit}{\end{itemize}}
\newcommand{\BNUM}{\begin{enumerate}}
\newcommand{\ENUM}{\end{enumerate}}
\newcommand{\BA}{\begin{array}}
\newcommand{\EA}{\end{array}}

\newcommand{\BSEQ}{\begin{subequations}}
\newcommand{\ESEQ}{\end{subequations}}
\newcommand{\BPM}{\begin{pmatrix}}
\newcommand{\EPM}{\end{pmatrix}}

\newcommand{\eg}{{\it e.g.}}
\newcommand{\ie}{{\it i.e.}}

\newcommand{\expect}{\mathop{\mathrm{E}}}

{\begin{quote}}{\end{quote}}

\begin{document}

\title{Critical behaviour in charging of electric vehicles}
\author{Rui Carvalho}
\email{rui.carvalho@durham.ac.uk}
\affiliation{School of Engineering and Computing Sciences, Durham University, Lower Mountjoy, South Road, Durham, DH1 3LE, UK}
\author{Lubos Buzna}
\email{lubos.buzna@fri.uniza.sk}
\affiliation{University of Zilina, Univerzitna 8215/1, 01026 Zilina, Slovakia}
\author{Richard Gibbens}
\email{richard.gibbens@cl.cam.ac.uk}
\affiliation{Computer Laboratory, University of Cambridge, William Gates Building, 15 JJ Thomson Avenue, Cambridge, CB3 0FD, UK}
\author{Frank Kelly}
\email{f.p.kelly@statslab.cam.ac.uk}
\affiliation{Statistical Laboratory, Centre for Mathematical Sciences, University of Cambridge, Wilberforce Road, Cambridge CB3 0WB, UK}
\begin{abstract}
The increasing penetration of electric vehicles over the coming decades, taken together with the high cost to upgrade local distribution networks and consumer demand for home charging, suggest that managing congestion on low voltage networks will be a crucial component of the electric vehicle revolution and the move away from fossil fuels in transportation. Here, we model the max-flow and proportional fairness protocols for the control of congestion caused by a fleet of vehicles charging on two real-world distribution networks. We show that the system undergoes a continuous phase transition to a congested state as a function of the rate of vehicles plugging to the network to charge. We focus on the order parameter and its fluctuations close to the phase transition, and show that the critical point depends on the choice of congestion protocol. Finally, we analyse the inequality in the charging times as the vehicle arrival rate increases, and show that charging times are considerably more equitable in proportional fairness than in max-flow. 
\end{abstract}
\maketitle 

\section{Introduction}


Electric vehicles may become competitive, in terms of total ownership costs, with internal-combustion engine vehicles over the next couple of decades. Studies in the United States and the UK suggest the current power grid has enough generation capacity to charge $70\%$ of cars and light trucks overnight, during periods of low demand~\cite{Service09}. A recent survey suggests, however, vehicle owners prefer home charging, would consider charging their vehicles during the day (typically between 6 and 10 pm), and are unwilling to accept a charging time of 8 hours~\cite{Deloitte10}. 
The time to fully charge the battery of an electric vehicle at home currently varies from 18 hours (Level 1, in the United States at 110 V and 15 A with a charge power of 1.4 kW) to 4 hours (Level 2, at 220 V, 30 A with a charge power of 6.6 kW). Alternatively, electric vehicles could charge at public charging stations at Level 3 in less than 30 minutes~\cite{Dickerman10}. Taken together, consumer behaviour and advances in battery technology may lead to a rise in  peak demand with the increasing penetration of electric vehicles, overloading distribution networks and requiring local infrastructure reinforcement~\cite{Clement10,Green11,Tran12,Keshav12}. To reduce the cost of upgrades to the last mile of cables, network operators may need to coordinate charging strategies in a way that is both technically and socially acceptable. 
To achieve this goal, network designers could implement charging protocols that prioritise the access of a fleet of electric vehicles to electric power, thus simultaneously managing network congestion and accounting for the fairness of user allocations.

Through a series of papers, the power grid has recently gained increased visibility in the scientific community~\cite{DSouza13,Baptista14}, and physicists have helped to increase our understanding of its synchronization~\cite{Motter13,Rohden12} and stability~\cite{Sole08,Kurths14}. In parallel, recent advances in optimization and phase transitions~\cite{Scala14,Seoane14}  suggest that the tools of critical phenomena and optimization can be merged, opening up new horizons. From the point of view of the distribution network operator, the problem of vehicle charging is to manage congestion on distribution networks, while respecting Kirchhoff's laws and keeping voltage drops bounded. Here, we explore two congestion control mechanisms: max-flow and proportional fairness. We show that if too many vehicles plug-in to the network, charging takes too long, more cars arrive than leave fully charged, and the system undergoes a continuous phase transition to a congested state~\cite{Guimera02,Lai05}, where the critical point depends on the choice of congestion control algorithm. By gaining insights into the critical behaviour that naturally emerges with the increase of the number of vehicles, we hope to help network designers decide which algorithms to implement in the real-world.

\section{The Model}
\label{sec:model}

Physicists are familiar with simulated annealing, a global optimization method that can avoid becoming trapped in a local optimum. In principle, it converges to the global optimum, but in practice this is not guaranteed (see \eg~\cite{Hajek88,Boyd01,Donetti05,Arenas08}) because the required theoretical cooling schedules are too slow to use in implementations. In contrast, convex optimization always finds the solution, if it exists, independently of the starting point. Convex optimization problems can be solved efficiently (typically in polynomial time), even for problems with hundreds of variables and thousands of constraints, using interior-point methods~\cite{Boyd04}. The burgeoning field of convex optimization in electricity networks~\cite{TaylorBook15,Low14b,Low14c} is a good example of an application of the mathematical framework developed over the last 20 years. Indeed, the extensive numerical simulations we present here are only possible due to techniques developed since 2012~\cite{Lavaei12,Low14b,Low14c}.
The networks that we study are relatively small. The stochasticity of vehicle arrival times, however, implies solving an optimization problem in each time step if the state of the system changes. Hence, to gain insights into the steady state of vehicle charging, efficient algorithms are a necessity at the design stage. Of course, real-world implementations also depend on efficient algorithms, which would need to run online, often in large urban distribution networks.  

An optimization problem is determined by a function of a set of variables (the objective function), for which we seek a minimum, and a set of upper bound constraints that restrict the domain (or \textit{feasible set}) of those variables~\cite{Boyd04}. A point is feasible if it belongs to the feasible set, and is optimal if it is the minimum of the objective function in the feasible set. 
An optimization problem is convex if both the objective function and the constraints are convex, in which case the objective function has a global minimum. A convex relaxation of an optimization problem $\mathsf{P}$ is a convex optimization problem $\mathsf{P}^\prime$ with an enlarged feasible set. If the optimum of $\mathsf{P}^\prime$ is feasible for $\mathsf{P}$, it is also the optimum for $\mathsf{P}$ and we say the relaxation is exact. Hence, convex relaxations are more attractive than approximate methods, such as linearisations, because the feasibility of the relaxed optimum of $\mathsf{P}^\prime$, which can be verified either analytically or numerically, is a certificate of the exactness of the relaxation.

Consider a tree topology, such that electric power is distributed from a root node to electric vehicles that charge at the nodes. Let $\mathcal{P}(t)$ be the feasible set of power allocations at time $t$, \ie~the set of all allocations of power to electric vehicles that do not violate the operational constraints of the distribution network. Each feasible allocation $P(t) \in \mathcal{P}(t)$ is a vector $P(t) = (P_l(t): l = 1, \dots, N(t))$, where $N(t)$ is the number of vehicles in the network at time $t$. Vehicle $l$ derives a utility $U_l(P_l(t))$ from the allocated charging power $P_l(t)$, and we wish to select the allocation that maximises the sum of vehicle utilities~\cite{Kelly14}. This allocation acts as a network protocol that distributes network capacity among users, and solves the following problem:
\BSEQ
\begin{empheq}{align}
\label{eq:max_flow_final_a}
\text{maximise} \quad & \sum\limits_{l=1}^{N(t)} U_l(P_l(t)) &&\\
\label{eq:max_flow_final_b}
\text{subject to } \quad & P(t) \in \mathcal{P}(t).
\end{empheq}
\ESEQ
Here we explore two user utility functions. First, we consider the non-unique \textit{max-flow} allocations given by $U_l\left( P_l(t) \right) = P_l(t)$, \ie~we maximise the instantaneous aggregate power sent from the root node to the vehicles, which is a benchmark of efficient network throughput~\cite{Bertsimas10}. Such allocations, however, can also leave users with zero power, which is considered unfair from the user point of view. Hence, we next consider the \textit{proportional fairness} allocation. 
Mathematically, the problem is to find the feasible allocation that maximises the sum of the logarithm of user rates, that is $U_l\left( P_l(t) \right)= log(P_l(t))$. The proportional fairness allocation is especial, because the users and the network operator simultaneously maximise their utility functions~\cite{Kelly14}. Furthermore, the problem is convex, and so can be solved in polynomial time~\cite{Boyd04}, and it can be naturally extended by adding positive weights to each term in the objective function Eq.~(\ref{eq:max_flow_final_a}), to account for diversity in user demand or for more than one user at some nodes~\cite{Kelly14}. 
For the compact and convex set $\mathcal{P}(t)$, it can be shown that the allocation $P^{PF}(t)$ that maximises Eq.~(\ref{eq:max_flow_final_a}), satisfies~\cite{Kelly14,Luss_2012}:
\BEQ
\label{eq:proportional_fairness_frac}
\sum_{l=1}^{N(t)}\frac{P_l(t) - P_l^{PF}(t)}{P_l^{PF}(t)} \leq 0.
\EEQ
This allocation is known as proportionally fair, because the aggregate of proportional changes with respect to all other feasible allocations is non-negative. In other words, Eq. (\ref{eq:proportional_fairness_frac}) implies that to increase the instantaneous power allocated to a vehicle by a percentage $\epsilon$, we have to decrease a set of other power allocations, such that the sum of the percentage decreases is larger or equal to $\epsilon$. In contrast, in max-flow, to increase the instantaneous power allocated to a vehicle by $\epsilon$, we have to decrease the sum of instantaneous powers allocated to other vehicles at least by $\epsilon$. 
It turns out that fairness and congestion control are two sides of the same coin, because the existing algorithms for fair allocations also manage network congestion~\cite{Bertsekas92,Kelly98,Tan99,Srikant04,Carvalho12,Carvalho14,Kelly14}.
Moreover, in the analysis of the parallel problem for communication networks, proportional fairness has emerged as a compromise between efficiency and fairness with an attractive interpretation in terms of shadow prices and a market clearing equilibrium~\cite[Section 7.2]{Kelly14}. 	

\begin{figure}
\centering
\includegraphics[width=\textwidth]{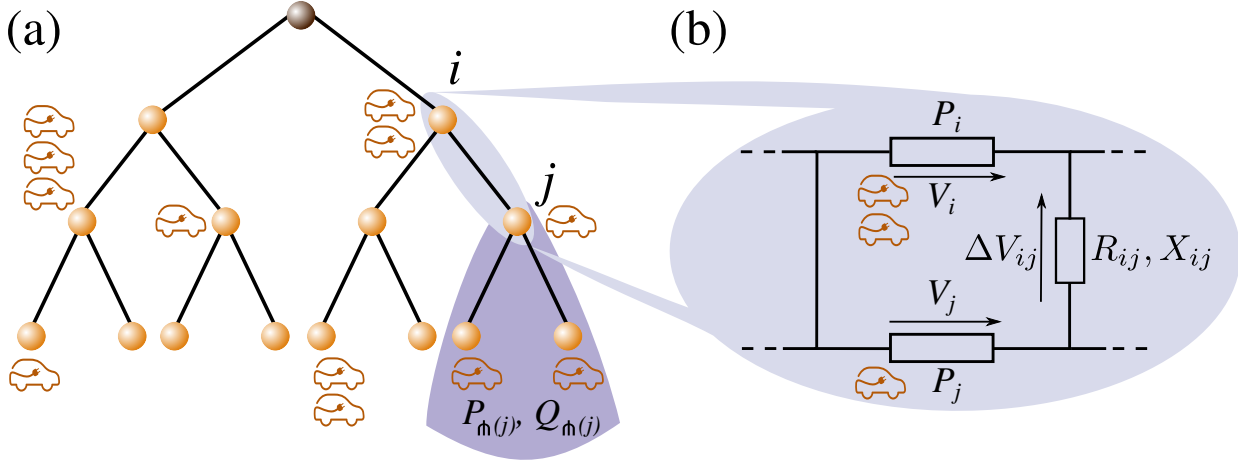}
\caption{Schematic illustration of (a) a distribution network, (b) the circuit of a network edge. Electric vehicles choose a charging node with uniform probability, and plug-in to the node until fully charged, as illustrated by the electric vehicle icons on the network. Network edge $e_{ij}$ has impedance $Z_{ij}=R_{ij}+iX_{ij}$. The power consumed by the subtree $\pitchfork(j)$ rooted at node $j$ (area shaded in purple) is $S_{\pitchfork(j)} = P_{\pitchfork(j)} + iQ_{\pitchfork(j)}$, where vehicles consume real power only, but network edges have both active (real) and reactive (imaginary) power losses.}
\label{Fig:tree_network}
\end{figure}

The simplest flow model in electricity networks is the DC power flow, which is a linearisation of the AC power flow equations, and thus can be solved with tools of linear programming. It assumes that voltage amplitude is constant on all nodes, a good approximation at the level of the high-voltage transmission network, but a poor one on local distribution networks. Indeed, voltage drops are significant in distribution networks, and determine the network capacity, which leads us to using models of power flow specific to distribution networks~\cite{Kersting01}.
We abstract the distribution network to a rooted directed tree $\pitchfork(r)$ with node (often called \textit{bus}) set $\mathcal{V}$, edge (also called \textit{branch}) set $\mathcal{E}$, and a root node $r$ (\textit{feeder}) that injects power into the tree~\footnote{We write $\pitchfork(r)$ instead of $\pitchfork(\mathcal{V},\mathcal{E},r)$ to simplify the notation.}. Edge $e_{ij}\in \mathcal{E}$ connects node $i$ to node $j$, where $i$ is closer to the root than $j$, and is characterised by the impedance $Z_{ij}=R_{ij}+iX_{ij}$, where $R_{ij}$ is the edge resistance and $X_{ij}$ the edge reactance. The power loss along edge $e_{ij}$ is given by $S_{ij}(t) = P_{ij}(t) + i Q_{ij}(t)$, where $P_{ij}(t)$ is the real power loss, and $Q_{ij}(t)$ the reactive power loss. Electric vehicle $l$ receives active power $P_l(t)$ until charged, but does not consume reactive power~\cite{Low11} ---see Fig.~\ref{Fig:tree_network}(a). 
The vector $V(t)$ denotes the voltage allocated to the nodes. 
The voltage drop $\Delta V_{ij}$ down the edge $e_{ij}$ is the difference between the amplitude of the voltage phasors $V_i$ and $V_j$, assuming node $i$ is closer to the root $r$ than node $j$~\cite{Kersting01}. Let $\pitchfork(j)$ denote the subtree of the distribution network rooted in node $j$, with node set $\mathcal{V}_{\pitchfork(j)}$ and edge set $\mathcal{E}_{\pitchfork(j)}$. Let $P_{\pitchfork (j)}$ denote the active power, and $Q_{\pitchfork(j)}$ the reactive power consumed by the subtree $\pitchfork(j)$---see Fig.~\ref{Fig:tree_network}.  Kirchhoff's voltage law applied to the circuit in Fig.~\ref{Fig:tree_network}(b) yields (see Appendix~\ref{appendix:voltage_drop}):
\BEQ
\label{eq:voltage_drop}
V_i(t)V_j(t) - V_j^2(t)  -P_{\pitchfork(j)}(t) R_{ij} - Q_{\pitchfork(j)}(t) X_{ij} = 0.
\EEQ

Vehicle $l$ has a battery with capacity $B$ that charges with the instantaneous power $P_l(t)$ from empty (at arrival time) to full (at departure time), and the level of battery charge is the time integral of instantaneous power. Vehicles arrive to the network, choose a node to charge randomly with uniform probability, charge until their battery is full, an lastly leave the network. At each time step, the network solves the congestion control problem to allocate instantaneous power to the vehicles.
The max-flow problem maximises the instantaneous aggregate power sent from the root node to the electric vehicles, respecting the constraints of distribution networks: the voltage drop along edges obeys Eq.~(\ref{eq:voltage_drop}), and node voltages are within $\left ( (1-\alpha)V_{nominal}, (1+\alpha)V_{nominal} \right )$ for $\alpha\in(0,1)$, with $\alpha=0.1$ typically~\cite{Kersting01}. Thus, to find the max-flow allocation of power to the vehicles, we solve the optimization problem for fixed $t$:
\BSEQ
\label{eq:problem_max-flow}
\begin{empheq}{align}
\mbox{$\underset{ V(t)}{\text{maximise}}$} \quad & U(t) =  \sum\limits_{l=1}^{N(t)} P_l(t) &&\label{eq:problem_max-flow_a}\\
\mbox{subject to} \quad & (1-\alpha) V_{nominal} \leq V_i(t) \leq (1+\alpha) V_{nominal}, && i \in \mathcal{V} \label{eq:problem_max-flow_b}\\
& V_i(t) V_j(t) -  V_j(t)V_j(t) - P_{\pitchfork (j)}(t) R_{ij} - Q_{\pitchfork (j)}(t) X_{ij} = 0 ,   && e_{ij} \in \mathcal{E}. \label{eq:problem_max-flow_c}
\end{empheq}
\ESEQ
Constraint~(\ref{eq:problem_max-flow_b}) sets the limits on the node voltage. Equation~(\ref{eq:problem_max-flow_c}) is the physical law coupling voltage to power, generalized from Eq.~(\ref{eq:voltage_drop}) for the subtree $\pitchfork (j)$, and encodes both Kirchhoff's voltage law on network edges and Kirchhoff's current law applied recursively at each node of the subtree (see Appendix~\ref{appendix:subtree}). We do not need to apply Kirchhoff's voltage law on network loops, however, because local distribution networks are approximately trees, and thus are cycle-free.  
Constraint~(\ref{eq:problem_max-flow_c}) is quadratic, hence not convex in general, which implies that the problem is not directly solvable by polynomial time methods. To overcome this limitation, we make a change of variables in problem~(\ref{eq:problem_max-flow}) by defining a weighted adjacency matrix $W(t)$, such that edge $e_{ij}$ corresponds to the $2\times 2$ principal submatrix $W(e_{ij},t)$ defined by~\cite{Low14a,BoseThesis14}:
\BEQ
\label{eq:def_W}
W(e_{ij},t) =  
\BPM
V_i(t)\\
V_j(t)
\EPM
\BPM
V_i(t) & V_j(t)
\EPM=
\BPM
V_i^2(t) & V_i(t)V_j(t) \\
V_j(t)V_i(t) & V_j^2(t) \\
\EPM=
\BPM
W_{ii}(t) & W_{ij}(t) \\
W_{ji}(t) & W_{jj}(t) \\
\EPM,
\EEQ
where $W_{ij}(t)=W_{ji}(t)$, because $V_i(t),V_j(t) \in \mathbb{R}$.
The matrices $W(e_{ij},t)$ are positive semidefinite, because their eigenvalues ($\lambda_1=0$ and $\lambda_2=V_i^2+V_j^2$) are non-negative, and rank one because they are of the form $v v^T$. Hence, constraint (\ref{eq:problem_max-flow_c}) can be replaced by three constraints: the first substitutes the quadratic terms in the voltages with linear terms in the $W(e_{ij},t)$, and the second and third constraints guarantee that the $W(e_{ij},t)$ are positive semidefinite and rank one. 

The solution of problem~(\ref{eq:problem_max-flow}) is on the \textit{Pareto frontier}~\footnote{We say that a power allocation $\{ P_l \}$ for $l=1,\ldots,N$ is better than another $\{ P^\prime_l \}$ if $P_l\ge P^\prime_l$ for all $l$, and for some $l$, $P_l>P^\prime_l$. A power allocation is \textit{Pareto optimal} or efficient if there is no better power allocation. The Pareto frontier of a set is the set of all Pareto optimal points.}, since we maximise an increasing function in the objective. 
The rank one constraint is nonconvex, but it does not change the Pareto frontier or the optimum~\cite{Lavaei14,Low14c}, and we remove it to relax problem~(\ref{eq:problem_max-flow}) to:
\BSEQ
\label{eq:max_flow_final}
\begin{empheq}{align}
\mbox{$\underset{W(t)}{\text{maximise}}$} \quad & U(t) =  \sum\limits_{l=1}^{N(t)} P_l(t) &&\\
\mbox{subject to} \quad & \left ( (1-\alpha) V_{\text{nominal}} \right )^2 \leq W_{ii}(t) \leq \left ( (1+\alpha) V_{\text{nominal}} \right )^2, && i \in \mathcal{V}\\
& W_{ij}(t) -  W_{jj}(t) - P_{\pitchfork (j)}(t) R_{ij} - Q_{\pitchfork (j)}(t) X_{ij} = 0,  && e_{ij} \in \mathcal{E}\\
& W(e_{ij},t) \succeq 0, && e_{ij} \in \mathcal{E}, && \label{eq:psd_W}
\end{empheq}
\ESEQ
where the generalized inequality in constraint~(\ref{eq:psd_W}) means the $W(e_{ij},t)$ matrices are positive semidefinite~\cite{Strang06}.

The problem of allocating power to vehicles in a proportional fair way has the same constraints as problem~(\ref{eq:max_flow_final}), however, the objective function is the sum of the logarithm of the power. It turns out, however, that it is computationally more efficient to aggregate vehicles at the nodes, and to maximise the sum of power allocated to the nodes, rather than the vehicles. To show this, we observe that all vehicles are equivalent, and thus the power $\mathrm{P_i}(t)$ allocated to node $i$ is divided equally among the vehicles charging on the node at each time step. Hence, if one or more vehicles is charging on node $i$, each gets the instantaneous power:
\BEQ
\label{eq:P_l}
P_l(t)=\frac{\mathrm{P_i}(t)}{w_i(t)},
\EEQ
where $w_i(t)= \sum_{l=1}^{N(t)} \Delta_{il}(t)$ is the number of electric vehicles charging on node $i$ at time $t$, and $\Delta_{il}(t)=1$ if electric vehicle $l$ is charging on node $i$ at time $t$ and zero otherwise. Hence, the proportional fair allocation is given by (see Appendix~\ref{appendix:aggregation_vehicles}):
\BSEQ
\label{eq:problem_proportional_fairness}
\begin{empheq}{align}
\mbox{$\underset{W(t)}{\text{maximise}}$} \quad & U(t)  = \sum_{i \in \mathcal{V}^+} w_i(t) \log \mathrm{P_i}(t) && \label{eq:problem_proportional_fairness_a}\\
\mbox{subject to} \quad & \left ( (1-\alpha) V_{\text{nominal}} \right )^2 \leq W_{ii}(t) \leq \left ( (1+\alpha) V_{\text{nominal}} \right )^2, && i \in \mathcal{V}\label{eq:problem_proportional_fairness_b}\\
& W_{ij}(t) -  W_{jj}(t) - P_{\pitchfork (j)}(t) R_{ij} - Q_{\pitchfork (j)}(t) X_{ij} = 0 , && e_{ij} \in \mathcal{E} \label{eq:problem_proportional_fairness_c}\\
& W(e_{ij},t) \succeq 0,  && e_{ij} \in \mathcal{E}, \label{eq:problem_proportional_fairness_d}
\end{empheq}
\ESEQ
where $\mathcal{V}^+$ is the subset  of  nodes with at least one vehicle, and we can recover the instantaneous power allocated to electric vehicle $l$, located at node $i$, from Eq.~(\ref{eq:P_l}). The complexity of the problem~(\ref{eq:problem_proportional_fairness}) thus scales with the number $\lvert \mathcal{V} \rvert$ of nodes, which is typically smaller than the number $N(t)$ of vehicles for large arrival rates $\lambda$. Similarly, we also aggregated vehicles in the implementation of problem~(\ref{eq:max_flow_final}), but omit the proof.

To study the behaviour of max-flow and proportional fairness as a function of the number of vehicles arriving at the network to be charged, we implement a discrete simulator that solves the congestion control problem in discrete time steps, starting with no vehicles charging on the network.
Vehicles arrive at the network in continuous time (following a Poisson process with rate $\lambda$) and with empty batteries, choose a node with uniform probability amongst all nodes (excluding the root), and charge at that node until their battery is full, at which point in time they leave the network. Once a vehicle plugs into a node, the congestion control algorithm will allocate it an instantaneous power, which is a function of the network topology and electrical elements, as well as the location of other vehicles. 

\begin{figure}
\centering
\subfigure{\includegraphics[width=0.495\textwidth]{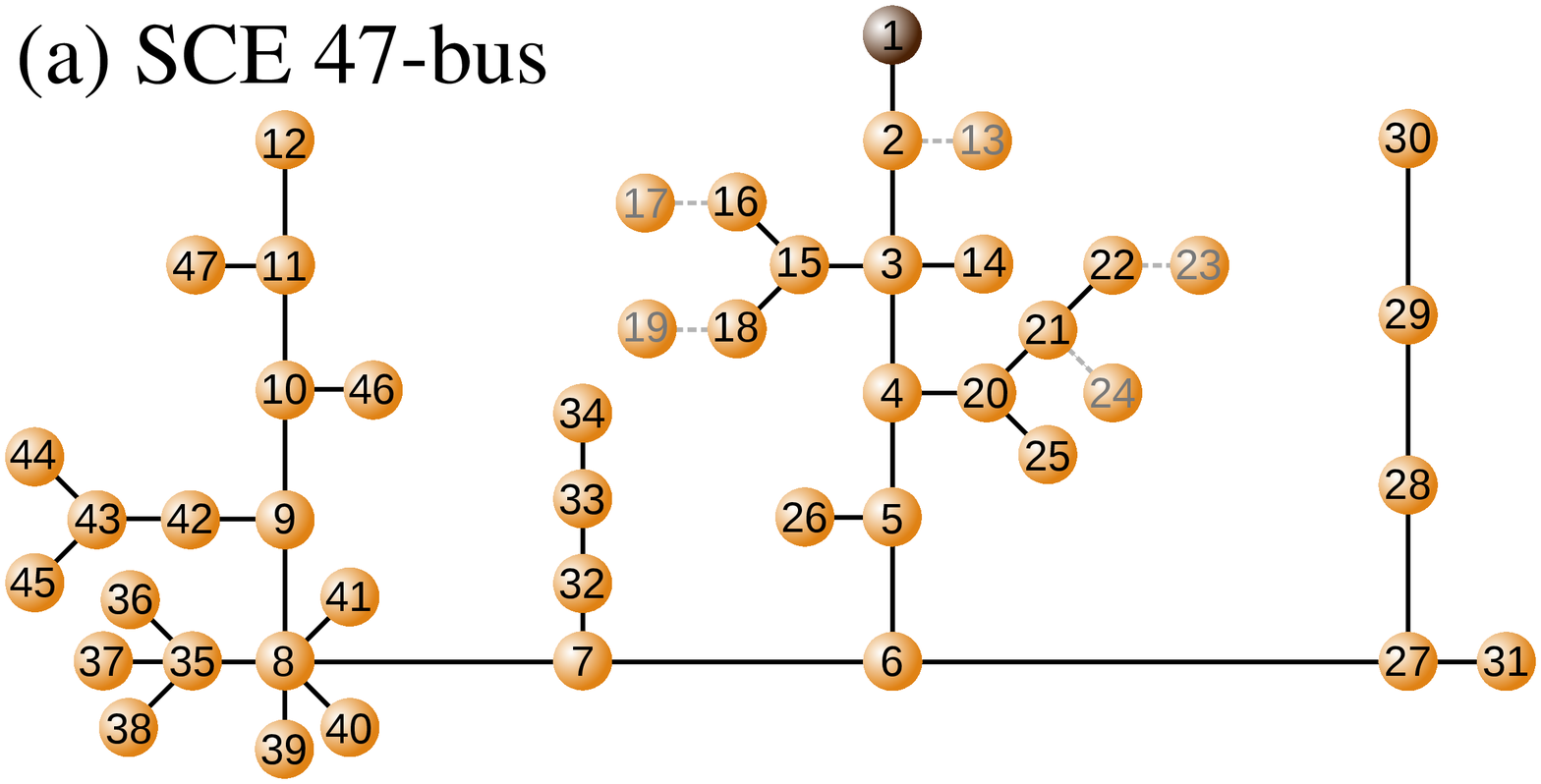}}
\subfigure{\includegraphics[width=0.495\textwidth]{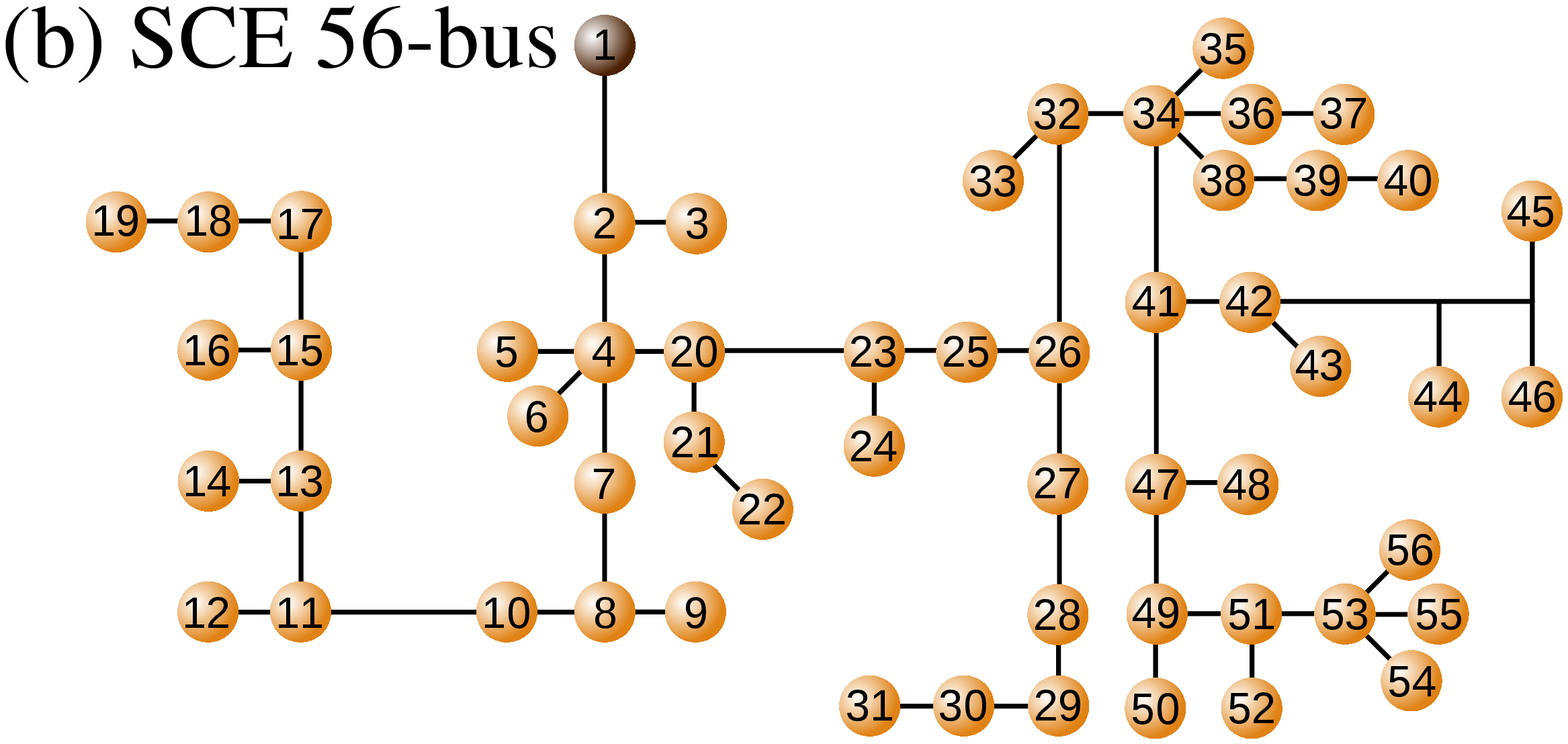}}
\caption{Topology of the (a) SCE 47-bus and (b) SCE 56-bus networks. Node indexes identify the edges, and edge resistance and reactance is taken from~\cite{Low14a}. Node 1 is the root node in both networks. Nodes 13, 17, 19, 23 and 24 of the SCE 47-bus network (in lighter colour) are photovoltaic generators, and we removed them from the network.}
\label{Fig:networks}
\end{figure}

 At each time step, we first check whether the number of charging vehicles changed (\ie~vehicles left the network fully charged, or new vehicles arrived to be charged), and if it has, we solve the max-flow problem~(\ref{eq:max_flow_final}) and the proportional fairness problem~(\ref{eq:problem_proportional_fairness}), which allocate a constant power during the time step to each of the charging vehicles. Next, we update the status of batteries at the end of the time step. The simulation terminates when the simulation time reaches the time horizon.
We simulated vehicles charging on the realistic SCE 47-bus and SCE 56-bus distribution networks~\cite{Low14a}, which are detailed in Fig.~\ref{Fig:networks}. To characterise the system behaviour in detail, we varied the arrival rate $\lambda$ from 0 to 1 in steps of 0.05 (0.005 close to the critical points), and for each $\lambda$ value we simulated an ensemble of 25 independent realisations of simulation runs, each simulation running for 15,000 time units (150,000 time units close to the critical point). We ran the simulations using CVXOPT~\cite{cvxopt} on the ETHZ Brutus cluster~\footnote{\url{https://www1.ethz.ch/id/services/list/comp_zentral/cluster/index_EN}} due to the high computational requirements. The computational time is comparable for max-flow and proportional fairness and for the 47-bus and 56-bus networks, but it is grows with $\lambda$.
For example, to simulate 5,000 time units of the proportional fairness algorithm for $\lambda = 1.0$ on the 47-bus network takes approximately 40 hours, but 4 minutes for $\lambda=0.05$ .

We set the battery capacity $B=1$ for all vehicles, and the nominal voltage $V_{nominal}=1$. Scaling $V_{nominal}$ by $\beta$, for $\beta\in (0,\infty)$, implies scaling both the the power delivered to vehicles and the battery capacity by $\beta^2$. To see this, observe that problems (\ref{eq:problem_max-flow}) and (\ref{eq:problem_proportional_fairness}) are invariant upon the scaling $V_{nominal}^\prime=\beta V_{nominal}$, $V_i^\prime(t)=\beta V_i(t)$ for all nodal voltages, $P_l^\prime(t)=\beta^2P_l(t)$,  $P_{\pitchfork}^\prime(t)=\beta^2P_{\pitchfork}(t)$ and $Q_{\pitchfork}^\prime(t)=\beta^2Q_{\pitchfork}(t)$, and $B^\prime=\beta^2 B$. 
Considering these scaling properties, our simulations can be extended to values of $V_{nominal}\neq 1$, provided the vehicle capacity $B$ is rescaled accordingly, and we use this property to rescale the problem when convenient.

\section{Numerical results}

\begin{figure}[h]
\centering
\subfigure{\includegraphics[width=\textwidth]{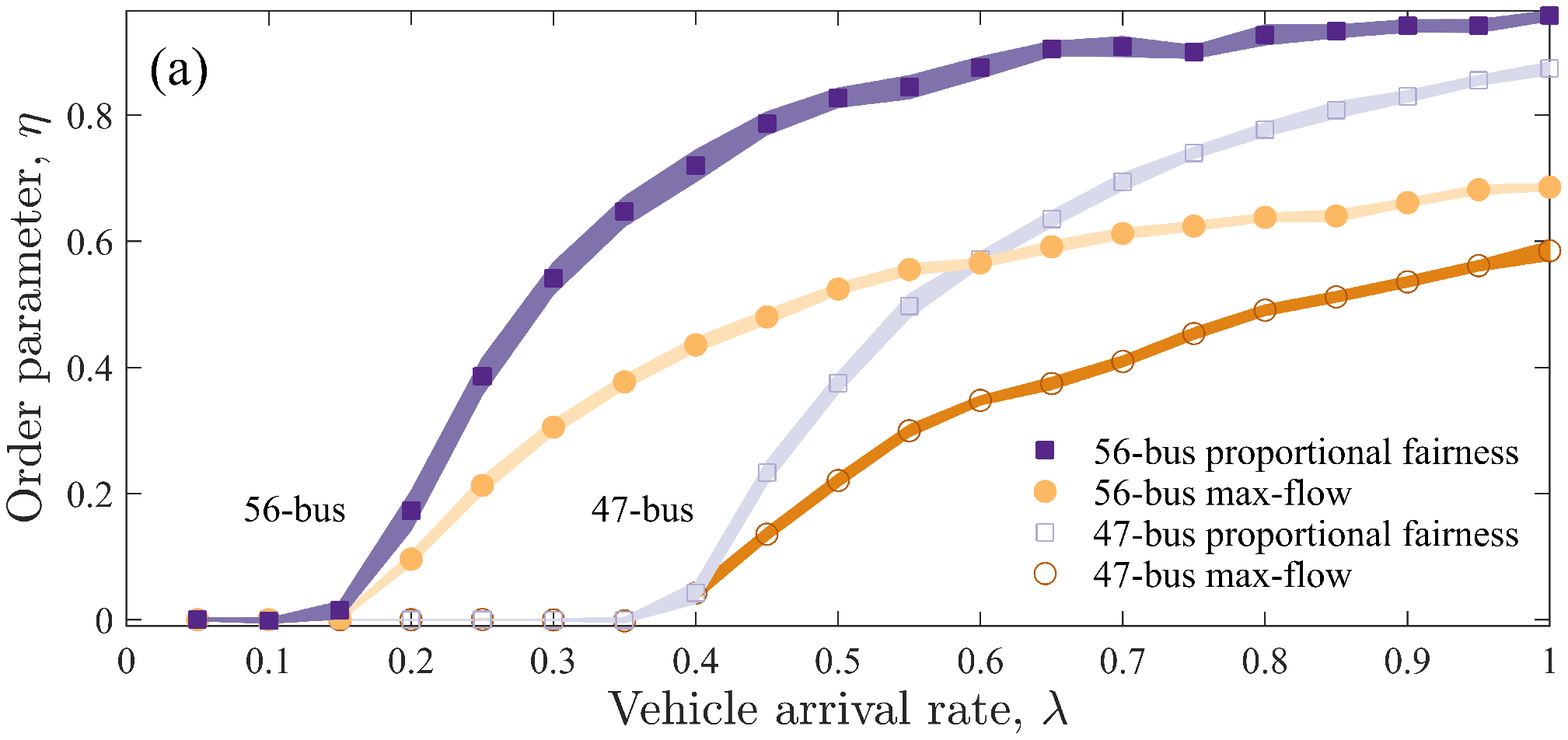}}\\
\subfigure{\includegraphics[width=0.45\textwidth]{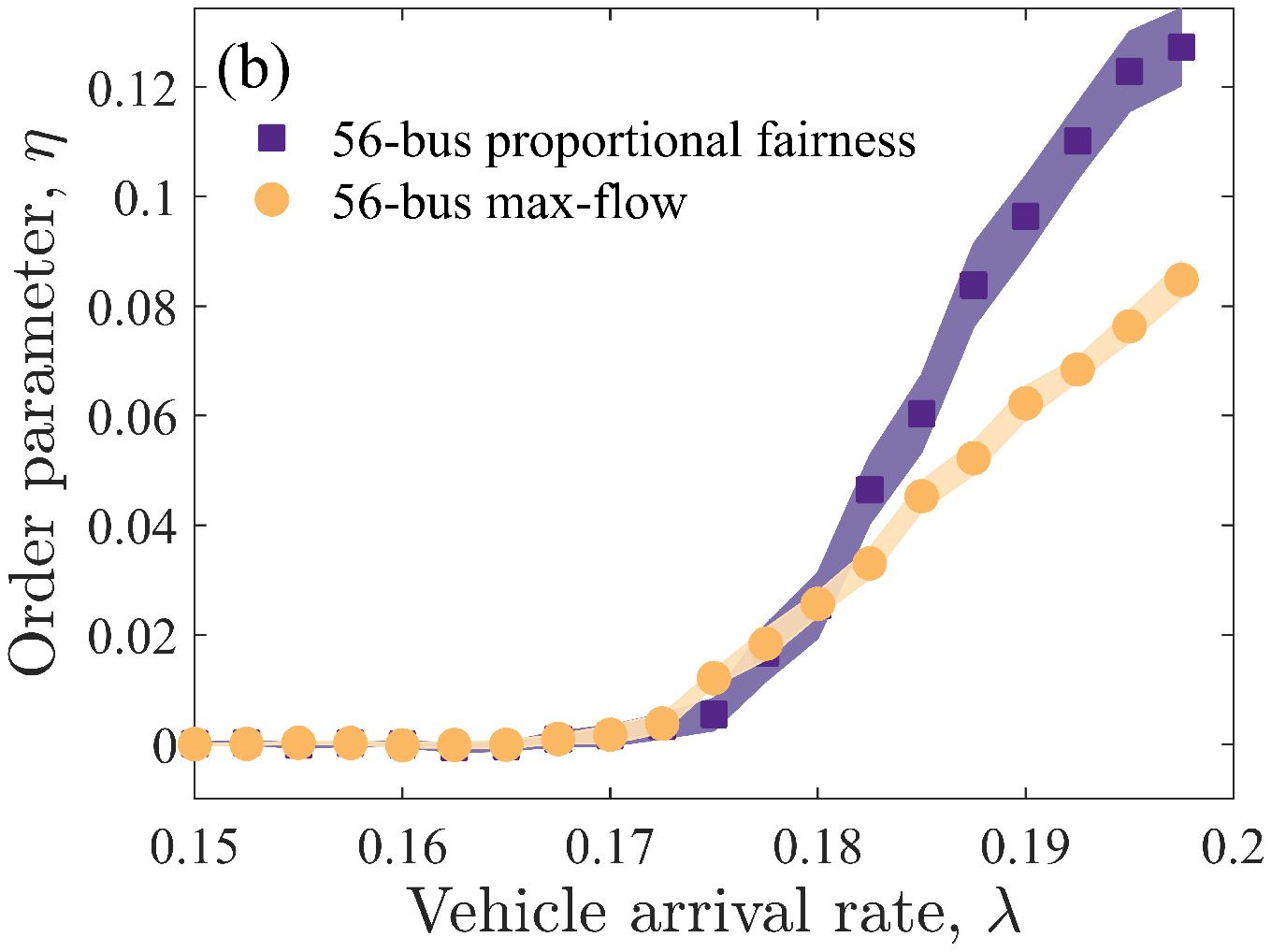}}
\subfigure{\includegraphics[width=0.45\textwidth]{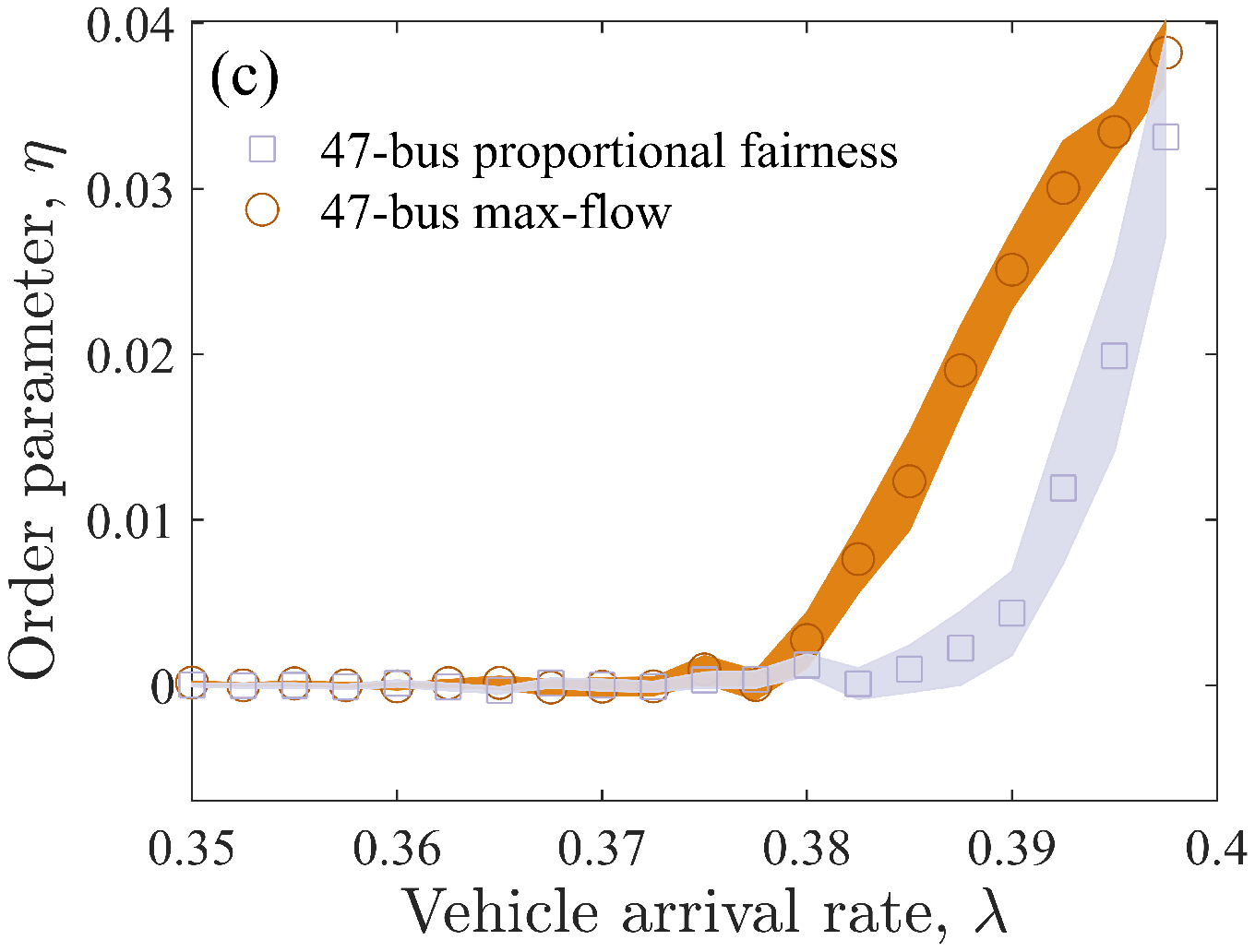}}
\caption{(a) Order parameter $\eta$ as a function of the arrival rate $\lambda$, for the SCE  56-bus (filled symbols) and 47-bus (unfilled symbols)  networks, where we apply the max-flow (circular symbols) and proportional fairness (square symbols) algorithms for the simulation horizon of $1.5\times 10^4$ time units. We plot a zoom of the critical region for the (b) 56-bus network and (c) 47-bus network for the longer horizon of $10^5$ time units.
Panel (c) suggests the critical arrival rate is different for the max-flow and proportional fairness algorithms in the 47-bus network. Symbols show average values over an ensemble of 25 runs and shaded areas represent $95\%$ confidence intervals.}
\label{fig:phase_diagram}
\end{figure}

We find critical behaviour that resembles results found in communication networks, in that both systems undergo a continuous phase transition~\cite{Guimera01}. In order to characterize this phase transition, we adopt the order parameter $\eta(\lambda)$ that represents the ratio at the steady state between the number of uncharged vehicles and the number of vehicles that arrive at the network to be charged~\cite{Guimera01}: 
\BEQ
\eta(\lambda) = \lim_{t\rightarrow\infty} \frac{1}{\lambda} \frac{\left\langle \Delta N(t)\right\rangle}{\Delta t},
\label{eq:eta}
\EEQ
where $\Delta N(t)=N(t+\Delta t) - N(t)$ and $\left\langle\ldots\right\rangle$ indicates an average over time windows of width $\Delta t$. We calculate $\eta(\lambda)$ in the steady state, that is $\lim_{t\rightarrow\infty} \Delta N(t)\propto  \Delta t$. For arrival rates $\lambda < \lambda_c$, all vehicles that plug-in to the network with empty batteries within a large enough time window leave fully charged within that period ({\it free flow} phase), but for $\lambda > \lambda_c$ some vehicles have to wait for increasingly long times to fully charge ({\it congested} phase). 
The order parameter characterises the phase transition:  $\eta(\lambda)=0$ in the free-flow regime, and $\eta(\lambda)>0$ in the congested phase, a higher order parameter meaning that queues of charging vehicles build up more rapidly.

Figure~\ref{fig:phase_diagram} is a plot of the order parameter for the 47 and 56-bus networks and the two congestion control methods, as a function of $\lambda$. Simulation results shown in Figs.~\ref{fig:phase_diagram}(a) suggest that $\lambda_c$ depends on several factors (the network topology, the complex impedance on the edges, battery capacity, $V_{nominal}$, as well as the position of vehicles on the network). At this resolution of the control parameter, it is unclear, however, whether the critical point is the same for max-flow and proportional fairness in both networks. To clarify this, we studied the order parameter with higher resolution close to the critical points---see Figs.~\ref{fig:phase_diagram}(b) and (c). The critical point is numerically indistinguishable for max-flow and proportional fairness in the 56-bus network. In the 47-bus network, however, we find that $\lambda_c$ is larger for proportional fairness than for max-flow.

\begin{figure}[h]
\centering
\subfigure{\includegraphics[width=0.45\textwidth]{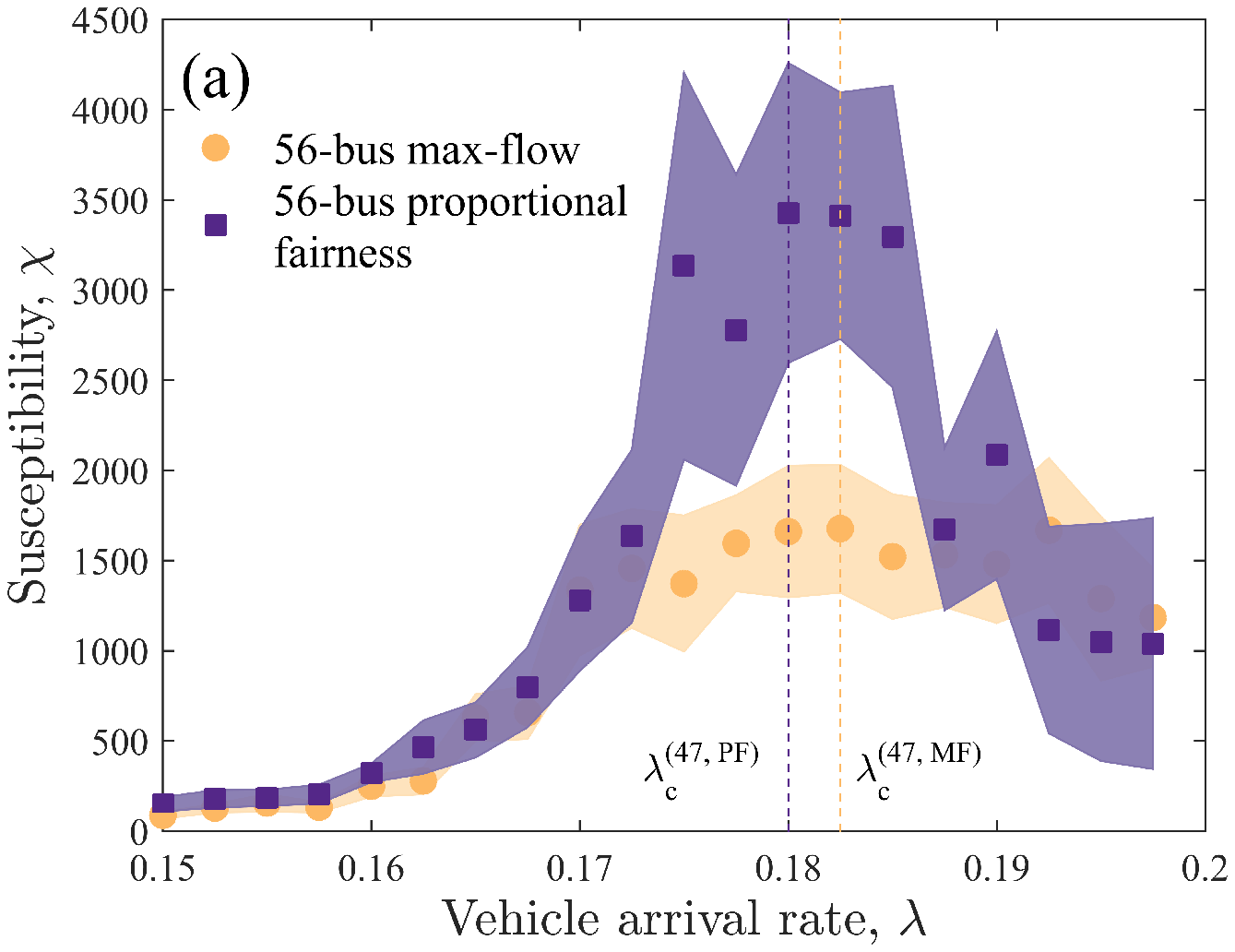}}
\subfigure{\includegraphics[width=0.45\textwidth]{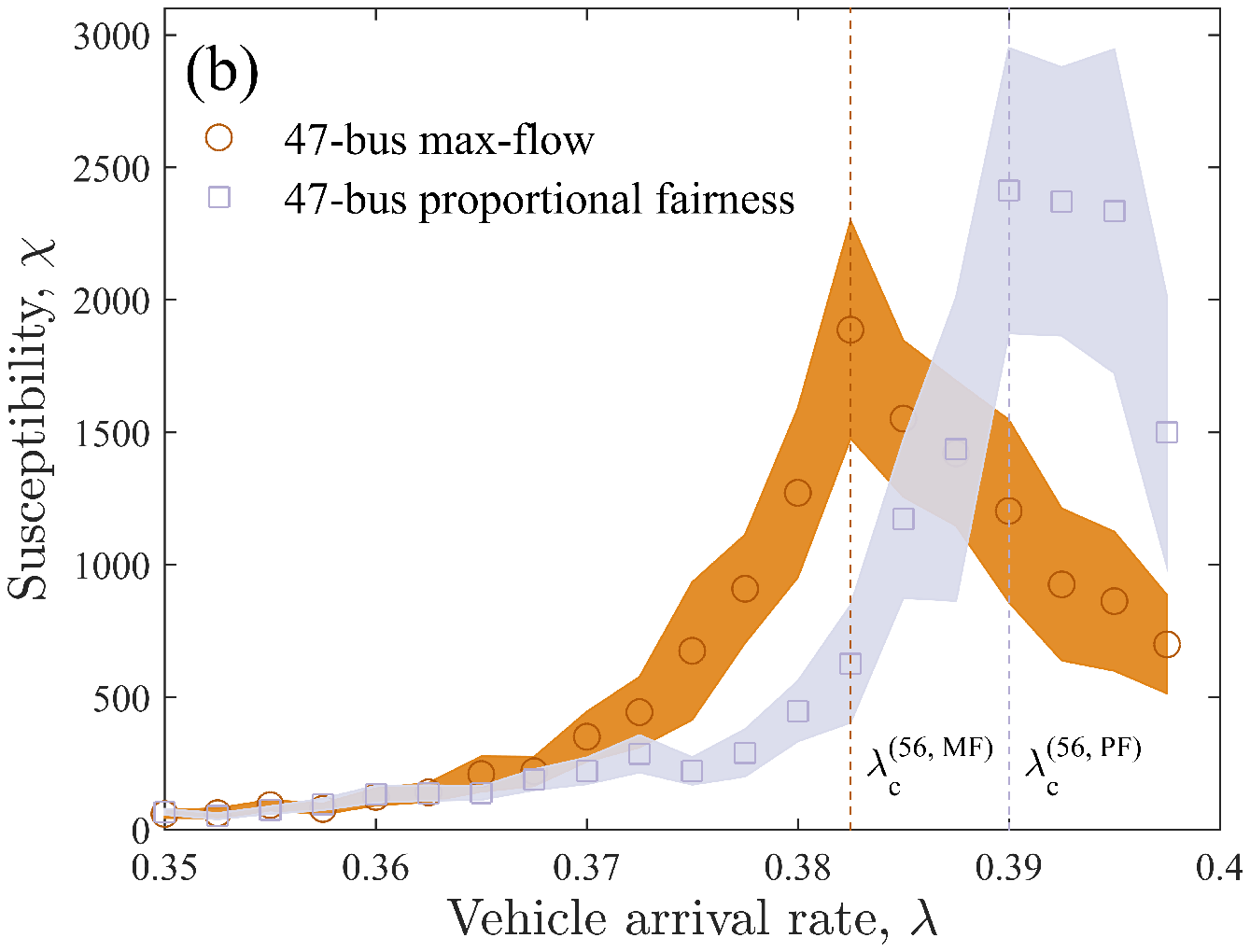}}
\subfigure{\includegraphics[width=0.45\textwidth]{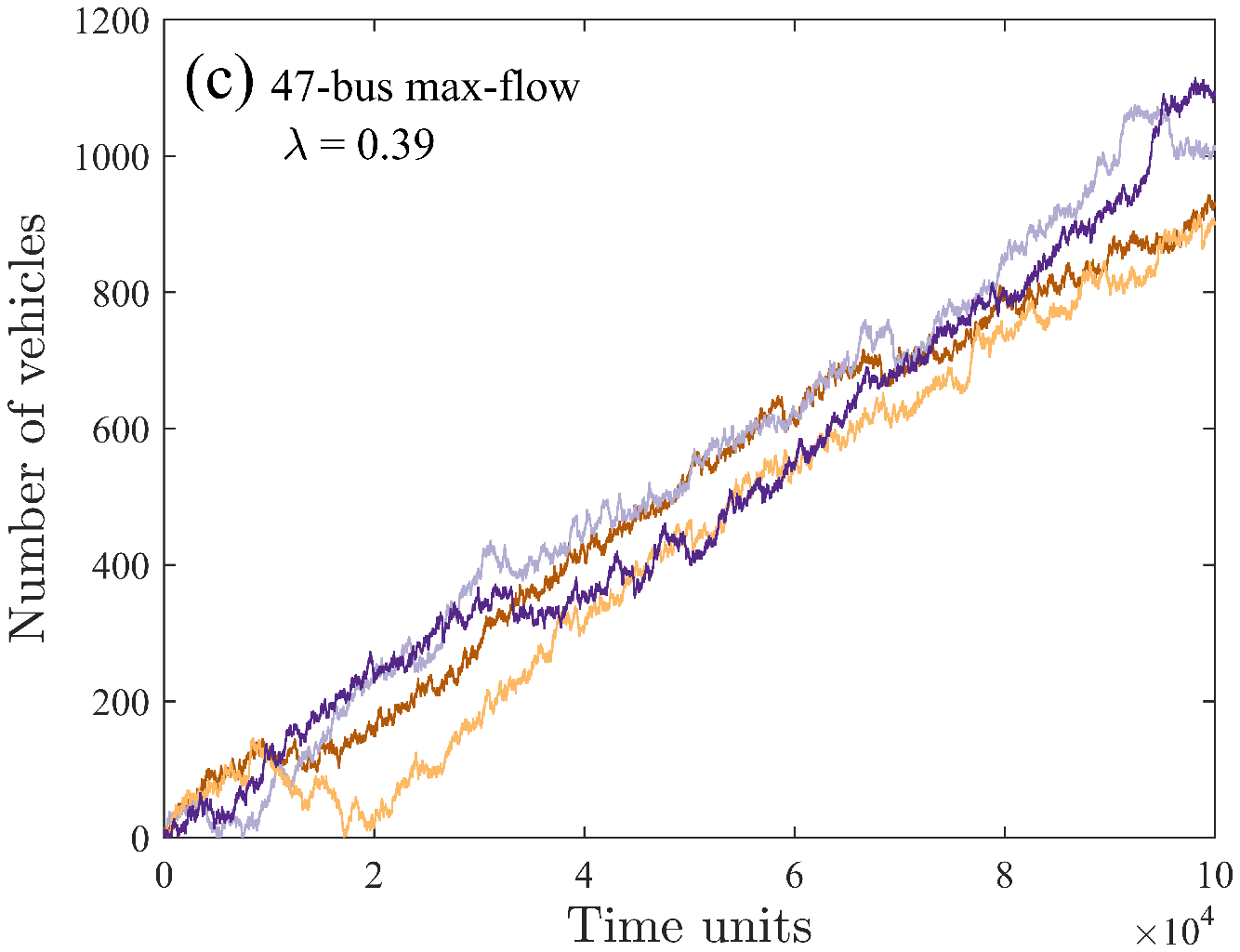}}
\subfigure{\includegraphics[width=0.45\textwidth]{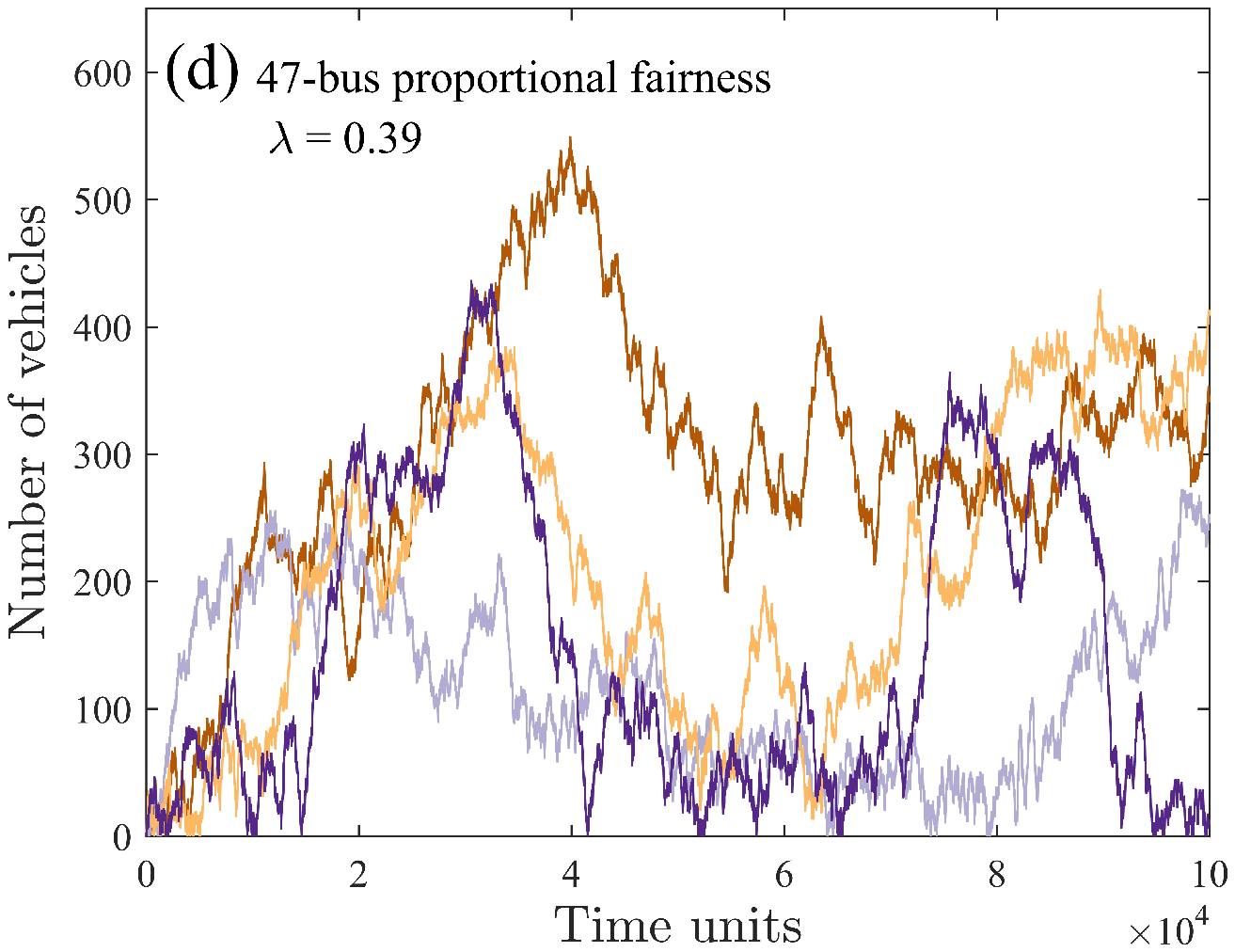}}
\caption{Susceptibility $\chi(\lambda)$ as a function of the arrival rate $\lambda$, for the for the (a) SCE  56-bus (filled symbols) and (b) 47-bus (unfilled symbols)  networks, where we apply the max-flow (circular symbols) and proportional fairness (square symbols) algorithms for the time horizon of $10^5$ time units. Vertical lines show the value of the critical points for max-flow (MF) and proportional fairness (PF). Panel (b) shows the difference in the critical arrival rate for the two congestion control algorithms. To illustrate this difference, we plot in (c) and (d) representative time series for the 47-bus network for $\lambda=0.39$, showing that, within the time horizon, max-flow is supercritical, whereas proportional fairness is subcritical. Symbols show average values over an ensemble of 25 runs and shaded areas represent $95\%$ confidence intervals.}
\label{fig:susceptibility}
\end{figure}

The number $N(t)$ of charging vehicles at time $t$ fluctuates widely close to the critical point, and thus it is difficult to determine $\lambda_c$ from Fig.~\ref{fig:phase_diagram}.
To overcome this limitation, we adopt the susceptibility-like function~\cite{Guimera01}:
\BEQ
\chi(\lambda)=\lim_{\Delta t \to \infty}\Delta t\ \sigma_\eta (\Delta t),
\EEQ
where $\Delta t$ is the length of a time window, and $\sigma_\eta (\Delta t)$ is the standard deviation of the order parameter $\eta$. To compute $\chi(\lambda)$, we first consider a long time series and split it into windows with length $\Delta t$. We next determine the value of the order parameter in each window, and finally calculate the standard deviation of these values. The susceptibility displays a singular point at $\lambda_c$ (see Fig.~\ref{fig:susceptibility})
, allowing us to study the dependencies of the critical arrival rates on the congestion control algorithm, as well as network topology and size. 

Similarly to our analysis of $\eta(\lambda)$, the values of $\lambda_c$ are indistinguishable in the 56-bus network. In contrast, however, in the 47-bus network the singular point of $\chi(\lambda)$ is smaller for max-flow than for proportional fairness. This suggests that proportional fairness charges a slightly larger number of vehicles than max-flow, and is thus marginally more efficient, on a neighbourhood of its critical point.
To support this conclusion, we show in Figs.~\ref{fig:susceptibility}(c) and (d) four representative instances of the time series of the number of vehicles charging on the 47-bus network at $\lambda=0.39$. The number $N(t)$ of vehicles grows linearly with time in max-flow in all four cases, suggesting that the critical point is below $\lambda=0.39$ for this algorithm. In contrast, $N(t)$ oscillates in proportional fairness, suggesting that the critical point is above $\lambda=0.39$, in agreement with the analysis of $\chi(\lambda)$. 

The two congestion control algorithms lead to different allocations of instantaneous power, with vehicles charging in different order and over different time intervals.  
If there are vehicles on a path $p$ between the root and a leaf node,  the voltage drops with increasing distance from the root, the lower limit voltage constraint (\ref{eq:problem_max-flow_b}) is fulfilled at equality for one node on $p$, and nodes further away than that will not receive any power. The objective function of proportional fairness guarantees that each vehicle gets a positive power allocation, thus the lower limit voltage constraint is satisfied at equality on the occupied node that is the most distant from the root on $p$. In max-flow, however, to maximise the aggregate power allocated to vehicles that can take all instantaneous power they are allocated (elastic demand), on a network with bounded voltage drops (\ie~capacity), implies also minimizing the power losses, and this is achieved by allocating all power on $p$ to the closest occupied node from the root on that path. For max-flow, this implies vehicles on the path $p$ further away from the root than the closest occupied node will only receive power after all vehicles on this node have left the network fully charged. 
In other words, under max-flow, users experience a charging time that depends strongly on their location on the network: vehicles close to the root charge faster, and vehicles on the tree leaves may take a very long time to charge. In contrast, under proportional fairness, the charging times are more homogeneous, because vehicles receive instantaneous powers that are also more uniform. 

\begin{figure}
\centering
\includegraphics[width=\textwidth]{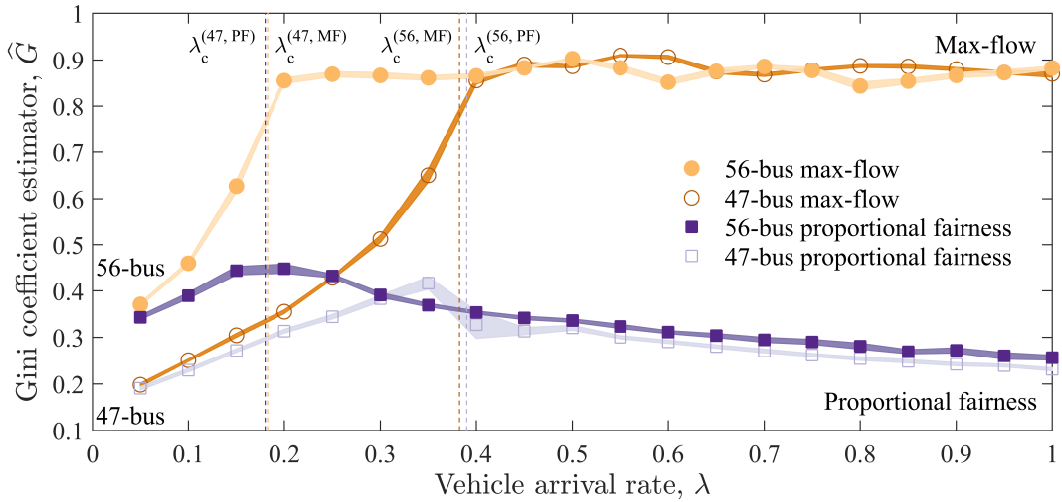}
\caption{Gini coefficient $G$ of the charging time as a function of the electric vehicle arrival rate $\lambda$, for the SCE 47-bus (unfilled symbols) and 56-bus (filled symbols) networks, where we apply the max-flow (circular symbols) and proportional fairness (square symbols) algorithms. We run the simulation for 15,000 times units, and compute the Gini coefficient from the charging time of vehicles that have charged fully during the simulation. To reduce the effect of a transient regime, we consider only vehicles that are fully charged after iteration $1000$. Vertical lines show the value of the critical points for max-flow (MF) and proportional fairness (PF) identified from the susceptibility $\chi(\lambda)$ for both networks. Symbols show average values over an ensemble of 25 runs and shaded areas represent $95\%$ confidence intervals.}
\label{Fig:Gini}
\end{figure}

To characterise inequalities in the user experience, we analyse the Gini coefficient of charging time. Originally devised as a measure of inequality in income distributions, the Gini coefficient is defined as~\cite{Ullah98}:
\BEQ
G=\frac{1}{2\mu} \expect \left [\left | u - v \right | \right ]=\frac{1}{2\mu}\int_{0}^{\infty}\int_{0}^{\infty}\left\vert u-v\right\vert f(u)f(v)\ du\ dv,
\label{eq:gini_coefficient_continuous}
\EEQ
where $u$ and $v$ are independent identically distributed random variables with probability density $f$ and mean $\mu$. In other words, the Gini coefficient is one half of the mean difference in units of the mean. The difference between the two variables receives a small weight in the tail of the distribution, where $f(u)f(v)$ is small, but a relatively large weight near the mode. Hence, $G$ is more sensitive to changes near the mode than to changes in the tails.
For a random sample ($x_i$, $i=1,2, \ldots, n$), the empirical Gini coefficient,~$\widehat{G}$, may be estimated by a sample mean
\BEQ
\widehat{G}=\frac{\sum_{i=1}^n\sum_{j=1}^n \left\vert x_i - x_j \right\vert}{2n^2\mu}.
\label{eq:gini_coefficient_discrete}
\EEQ
The Gini coefficient is used as a measure of inequality, because a sample where the only non-zero value is $x$ has $\mu=x/n$ and hence $\widehat{G}=(n-1)/n \to 1$ as~$n\to\infty$, whereas $\widehat{G}=0$ when all data points have the same value. 

We observe in Fig.~\ref{Fig:Gini} that the Gini coefficient of the charging time is larger in max-flow than in proportional fairness, for each of the networks. 
Moreover, the Gini coefficient increases faster in max-flow than in proportional fairness in the non-congested regime, showing that, when the system is stable, vehicles will experience a faster increase in the inequality of charging times in max-flow than in proportional fairness, with the increase of the vehicle arrival rate $\lambda$. 
For comparison with well-known measures of income inequality, Sweden has a Gini of 0.26, the United States has a Gini of 0.41 and the Seychelles has the highest Gini of 0.66~\cite{WorldBank15}. The proportional fairness algorithm reaches a maximum Gini of 0.45, which is comparable with the level of inequality in the US society, and thus may be judged sociable acceptable. The max-flow algorithm, however, reaches a Gini of 0.91, which measures a level of inequality considerably higher than present in any contemporary society. 

\section{Discussion}
In conclusion, we modelled the max-flow and proportional fairness protocols for the control of congestion caused by a fleet of vehicles charging on distribution networks. We analysed the second order phase transition that occurs with the increase of the number of electric vehicles  that arrive at the network with empty batteries to be charged, and found that the critical arrival rate $\lambda_c$ depends on the congestion control method. Indeed, we showed numerically on the 47-bus bus network that the onset of congestion takes place for larger values of $\lambda$ in proportional fairness than in max-flow.
This result is surprising, because one would expect that, for a chosen arrival rate $\lambda$, the maximisation of the aggregate instantaneous power would also lead to a maximisation of the energy (the time integral of power), and hence to a maximisation of the number of charged vehicles. This discovery, illustrates how the greediness of max-flow can be sub-optimal in relation to proportional fairness, which is an example of a fair allocation of instantaneous power.

We analysed the inequality in the charging times as the vehicle arrival rate increases, and showed that charging times are considerably more equitable in proportional fairness than in max-flow. Indeed, vehicles close to the root get all the power allocation in max-flow, leaving other vehicles excluded from the network and unable to charge. Hence, proportional fairness is preferable to max-flow, not only because it does not exclude users from the network, but also because the charging times are more equitable, and it can serve a higher number of vehicles. In conclusion, proportional fairness is a promising candidate protocol to manage congestion in the charging of electric vehicles.

\acknowledgements
We thank Janusz Bialek, Chris Dent and Emmanouil Loukarakis for help with modelling distribution networks. We also thank Dirk Helbing for granting access to the ETHZ Brutus high-performance cluster. This work was supported by the Engineering and Physical Sciences Research Council under grant number EP/I016023/1, by APVV (project APVV-0760-11) and by VEGA (project 1/0339/13).

\appendix
\section{Voltage drop on one edge}
\label{appendix:voltage_drop}
\begin{figure}[h]
\centering
\includegraphics[width=\textwidth]{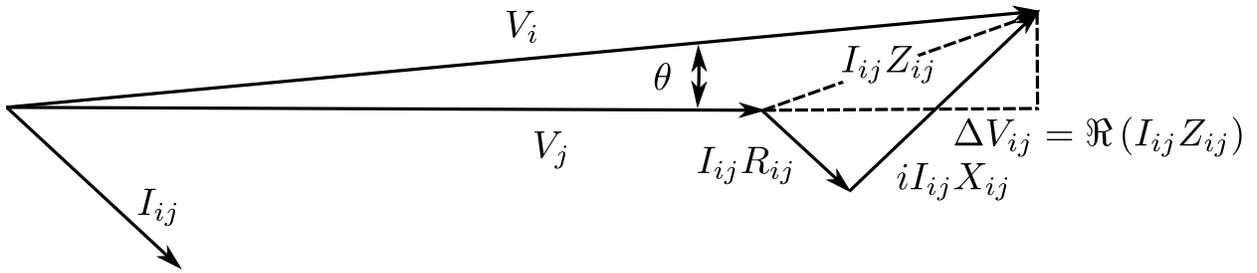}
\caption{The difference $I_{ij} Z_{ij}$ between the $V_i$ and $V_j$ phasors, decomposed  along the $V_j$ vector and its orthogonal direction. The phase angle $\theta$ difference between $V_i$ and $V_j$ is small, and hence the voltage drop can be approximated by $\Delta V_{ij}\simeq \Re\left(I_{ij}\text{ }Z_{ij}\right)$.} 
\label{Fig:phasor}
\end{figure}

The angle $\theta$ between $V_i$ and $V_j$ is small in distribution networks~\cite{Kersting01} (see Fig.~\ref{Fig:phasor}), and hence the phases of $V_i$ and $V_j$ are approximately the same, and can be chosen so the phasors have zero imaginary components.
Since the phasors are real, we can derive the voltage drop from Kirchhoff's voltage law applied to the circuit in Fig.~\ref{Fig:tree_network}(b),
\begin{align}
\label{eq:voltage_dropone}
 \Delta V_{ij} &= \rvert V_i\rvert- \rvert V_j\rvert \simeq V_i - V_j \nonumber=\\
& =\Re\left(I_{ij}\text{ }Z_{ij}\right) = \Re\left(\frac{I_{ij}\text{ }Z_{ij}\text{ }V_j^{*}}{V_j^{*}}\right)  =\Re\left(\frac{S_{\pitchfork(j)}^*\text{ }Z_{ij}} {V_j^{*}} \right)= \nonumber\\
& = \Re\left(\frac{(P_{\pitchfork(j)} - iQ_{\pitchfork(j)})(R_{ij} + iX_{ij})}{V_j^{*}}\right)  \simeq \frac{P_{\pitchfork(j)} R_{ij} + Q_{\pitchfork(j)} X_{ij} }{V_j},
\end{align}
where the superscript asterisk denotes the complex conjugate transpose.

\section{Active and reactive loads on a subtree}
\label{appendix:subtree}
From Kirchhoff's current law, the active and reactive power consumed by the loads in the subtree rooted in node $k$ can be computed as:
\begin{align}
\label{P_tree}
& P_{\pitchfork (k)} = \sum_{i \in \mathcal{V}_{\pitchfork(k)}} \sum_{l=1}^{N(t)} \Delta_{il}(t) P_l(t) + \sum_{i \in \mathcal{V}_{\pitchfork(k)}} \sum_{j: e_{ij}\in \mathcal{E}_{\pitchfork(k)}} P_{ij}(t),
\end{align}
and
\begin{align}
\label{Q_tree}
& Q_{\pitchfork (k)} =  \sum_{i \in \mathcal{V}_{\pitchfork(k)}} \sum_{j: e_{ij}\in \mathcal{E}_{\pitchfork(k)}} Q_{ij}(t),
\end{align}
where $P_{ij}(t)$ is the active and $Q_{ij}(t)$ the reactive power dissipated on a cable connecting nodes $i$ and $j$. The complex power is given by:
\begin{align}
& S_{ij}(t) = P_{ij}(t) + iQ_{ij}(t) = (V_i(t) - V_j(t)) I^{*} =  \nonumber \\ 
& = (V_i(t) - V_j(t)) \left( \frac{V_i(t) - V_j(t)}{R_{ij} + iX_{ij}} \right)^{*} = \nonumber \\
& = \left (V_i(t) - V_j(t)\right ) \left (V_i(t) - V_j(t)\right )^{*}\left( \frac{R_{ij} - iX_{ij}}{R_{ij}^{2} + X_{ij}^{2}} \right)^{*} = \nonumber \\
& = \left (W_{ii}(t)-W_{ij}(t)-W_{ji}(t)+W_{jj}(t)\right ) \frac{R_{ij} + iX_{ij}}{R_{ij}^{2} + X_{ij}^{2}}.
\end{align}
Since, the voltages are real, $W_{ij}(t)=W_{ji}(t)$, and thus
\BEQ
 P_{ij}(t) = \left (W_{ii}(t) -2W_{ij}(t) + W_{jj}(t)\right  )\frac{R_{ij}}{R_{ij}^{2} + X_{ij}^{2}},
\EEQ
and
\BEQ
Q_{ij}(t) = \left (W_{ii}(t) -2W_{ij}(t) + W_{jj}(t)\right  )\frac{X_{ij}}{R_{ij}^{2} + X_{ij}^{2}}. 
\EEQ

\section{Aggregation of vehicles at the nodes}
\label{appendix:aggregation_vehicles}
In proportional fairness, we maximise the sum of the logarithm of the instantaneous power allocated to electric vehicles:
\BEQ
\label{eq:utility_1}
\sum\limits_{l=1}^{N(t)}  \log P_l(t) = \sum_{i \in \mathcal{V^+}} \sum_{l=1}^{N(t)} \Delta_{il}(t) \log \frac{\mathrm{P_i}(t)}{\sum_{l=1}^{N(t)} \Delta_{il}(t)},
\EEQ
where $P_l(t)$ is the instantaneous power allocated to electric vehicle $l$, and $\mathrm{P_i}$ the instantaneous power allocated to node $i$. 
To maximise Eq.~(\ref{eq:utility_1}), we solve a problem with gradient and Hessian matrices that grow in size with the number of electric vehicles on the network. A more efficient way to approach the problem is to aggregate cars for each node $i$, then solve the optimization problem for the nodes (as if they were `super-cars'), and finally distribute the power allocated to each node among the cars on the node. To do this, we remove constant terms in the objective function Eq.~(\ref{eq:utility_1}), yielding:
\BEQ
\label{eq:utility_2}
U(t)  =  \sum_{i \in \mathcal{V}^+} \sum_{l=1}^{N(t)} \Delta_{il}(t) \log \mathrm{P_i}(t) = \sum_{i \in \mathcal{V}^+} w_i(t) \log \mathrm{P_i}(t).
\EEQ

\providecommand{\newblock}{}

\end{document}